\documentclass[11pt]{article}

\usepackage{theorem}

\newtheorem{theorem}{Theorem}

\newtheorem{corollary}[theorem]{Corollary}

{\theorembodyfont{\rmfamily}
 
 \newtheorem{example}{Example}
 \newtheorem{algorithm}{Algorithm}
}

\def\qed{\hfill $\Box$}

\newenvironment{proof}{\noindent{\bfseries Proof.}}{\qed}

\usepackage{amsmath}
\usepackage{amssymb}
\usepackage{latexsym}
\usepackage{amsfonts}

\title{\textbf{Words restricted by patterns with at most 2 distinct letters}}

\author{
{\bf Alexander Burstein}\\
{\small Department of Mathematics}\\
{\small Iowa State University}\\
{\small Ames, IA 50011-2064 USA}\\
{\small \texttt{burstein@math.iastate.edu}}\\
\and
{\bf Toufik Mansour}\\
{\small LaBRI, Universit\'e Bordeaux}\\
{\small 351 cours de la Lib\'eration}\\
{\small 33405 Talence Cedex, France}\\
{\small \texttt{toufik@labri.fr}} }

\begin{document}

\date{\today}

\maketitle

\begin{abstract}
We find generating functions for the number of words avoiding
certain patterns or sets of patterns with at most 2 distinct
letters and determine which of them are equally avoided. We also
find exact number of words avoiding certain patterns and provide
bijective proofs for the resulting formulas.
\end{abstract}


Let $[k]=\{1,2,\dots,k\}$ be a (totally ordered) alphabet on $k$
letters. We call the elements of $[k]^n$ \emph{words}. Consider
two words, $\sigma\in[k]^n$ and $\tau\in[\ell]^m$. In other words,
$\sigma$ is an $n$-long $k$-ary word and $\tau$ is an $m$-long
$\ell$-ary word. Assume additionally that $\tau$ contains all
letters $1$ through $\ell$. We say that $\sigma$ contains an
\emph{occurrence} of $\tau$, or simply that $\sigma$
\emph{contains} $\tau$, if $\sigma$ has a subsequence
\emph{order-isomorphic} to $\tau$, i.e. if there exist $1\le
i_1<\ldots<i_m\le n$ such that, for any relation
$\phi\in\{<,=,>\}$ and indices $1\le a,b \le m$,
$\sigma({i_a})\phi\sigma({i_b})$ if and only if $\tau(a) \phi
\tau(b)$. In this situation, the word $\tau$ is called a
\emph{pattern}. If $\sigma$ contains no occurrences of $\tau$, we
say that $\sigma$ \emph{avoids} $\tau$.

Up to now, most research on forbidden patterns dealt with cases
where both $\sigma$ and $\tau$ are permutations, i.e. have no
repeated letters. Some papers (Atkinson et al. \cite{Atkinson},
Burstein \cite{Burstein}, Regev \cite{Regev}) also dealt with
cases where only $\tau$ is a permutation. In this paper, we
consider some cases where forbidden patterns $\tau$ contain
repeated letters. Just like \cite{Burstein}, this paper is
structured in the manner of Simion and Schmidt \cite{SS}, which
was the first to systematically investigate forbidden patterns and
sets of patterns.


\section{Preliminaries}

Let $[k]^n(\tau)$ denote the set of $n$-long $k$-ary words which
avoid pattern $\tau$. If $T$ is a set of patterns, let $[k]^n(T)$
denote the set of $n$-long $k$-ary words which simultaneously
avoid all patterns in $T${\!}, that is $[k]^n(T)=\cap_{\tau\in
T}[k]^n(\tau)$.

For a given set of patterns T, let $f_T(n,k)$ be the number of
$T$-avoiding words in $[k]^n$, i.e. $f_T(n,k)=|[k]^n(T)|$. We
denote the corresponding exponential generating function by
$F_T(x;k)$; that is, $F_T(x;k)=\sum_{n\geq 0}{f_T(n,k)x^n/n!}$.
Further, we denote the ordinary generating function of $F_T(x;k)$
by $F_T(x,y)$; that is, $F_T(x,y)=\sum_{k\geq 0} F_T(x;k)y^k$. The
reason for our choices of generating functions is that $k^n\ge
|[k]^n(T)|\ge n!\binom{k}{n}$ for any set of patterns with
repeated letters (since permutations without repeated letters
avoid all such patterns). We also let
$G_T(n;y)=\sum_{k=0}^{\infty}{f_T(n,k)y^k}$, then $F_T(x,y)$ is
the exponential generating function of $G_T(n;y)$.

We say that two sets of patterns $T_1$ and $T_2$ belong to the same
\emph{cardinality class}, or \emph{Wilf class}, or are
\emph{Wilf-equivalent}, if for all values of $k$ and $n$, we have
$f_{T_1}(n,k)=f_{T_2}(n,k)$.

It is easy to see that, for each $\tau$, two maps give us patterns
Wilf-equivalent to $\tau$. One map, $r:\tau(i)\mapsto\tau(m+1-i)$,
where $\tau$ is read right-to-left, is called \emph{reversal}; the
other map, where $\tau$ is read upside down,
$c:\tau(i)\mapsto\ell+1-\tau(i)$, is called
\emph{complementation}. For example, if $\ell=3$, $m=4$, then
$r(1231)=1321$, $c(1231)=3213$, $r(c(1231))=c(r(1231))=3123$.
Clearly, $c\circ r=r\circ c$ and $r^2=c^2=(c\circ r)^2=id$, so
$\langle r,c\rangle$ is a group of symmetries of a rectangle.
Therefore, we call $\{\tau,\, r(\tau),\, c(\tau),\, r(c(\tau))\}$
the \emph{symmetry class} of $\tau$.

Hence, to determine cardinality classes of patterns it is enough
to consider only representatives of each symmetry class.


\section{Two-letter patterns}

There are two symmetry classes here with representatives 11 and 12.
Avoiding 11 simply means having no repeated
 letters, so
\[ \label{f11}
f_{11}(n,k)=\binom{k}{n}n!=(k)_n,\qquad
F_{11}(x;k)=\left(1+x\right)^k.
\]
A word
avoiding 12 is just a non-increasing string, so
\[ \label{f12}
f_{12}(n,k)=\binom{n+k-1}{n},\qquad
F_{12}(x;k)=\frac{1}{(1-x)^k}.
\]


\section{Single 3-letter patterns}

The symmetry class representatives are 123, 132, 112, 121, 111. It
is well-known \cite{Knuth} that
\[
|S_n(123)|=|S_n(132)|=C_n=\frac{1}{n+1}\binom{2n}{n},
\]
the $n$th Catalan number. It was also shown earlier by the first
author \cite{Burstein} that
\[
f_{123}(n,k)=f_{132}(n,k)=2^{n-2(k-2)}\sum_{j=0}^{k-2}{a_{k-2,j}\binom{n+2j}{n}},
\]
where
\[
a_{k,j}=\sum_{m=j}^{k}{C_m D_{k-m}}, \quad D_t=\binom{2t}{t},
\]
and
\[
F_{123}(x,y)=F_{132}(x,y)=1+\frac{y}{1-x}+
\frac{2y^2}{(1-2x)(1-y)+\sqrt{((1-2x)^2-y)(1-y)}} \, .
\]

\bigskip

Avoiding pattern 111 means having no more than 2 copies of each
letter. There are $0\le i\le k$ distinct letters in each word
$\sigma\in[k]^n$ avoiding 111, $0\le j\le i$ of which occur twice.
Hence, $2j+(i-j)=n$, so $j=n-i$. Therefore,
\[ \label{f111b}
f_{111}(n,k)=\sum_{i=0}^{k}
{\binom{k}{i}\binom{i}{n-i}\frac{n!}{2^{n-i}}}
=\sum_{i=0}^{k}{\frac{n!}{2^{n-i}(n-i)!(2i-n)!}(k)_i}=\sum_{i=0}^{k}{B(i,n-i)(k)_i},
\]
where $(k)_i$ is the falling factorial, and $B(r,s)=\displaystyle
\frac{(r+s)!}{2^s (r-s)!s!}$ is the Bessel number of the first
kind. In particular, we note that $f_{111}(n,k)=0$ when $n>2k$.
\begin{theorem}
\label{f111}
\(
\displaystyle F_{111}(x;k)=\left(1+x+\frac{x^2}{2} \right)^k.
\)
\end{theorem}
\begin{proof}
This can be derived from the exact formula above.
Alternatively, let $\alpha$ be any word in $[k]^n(111)$. Since
$\alpha$ avoids $111$, the number of occurrences of the letter $k$
in $\alpha$ is 0, 1 or 2. Hence, there are $f_{111}(n,k-1)$,
$nf_{111}(n-1,k-1)$ and $\binom{n}{2}f_{111}(n-2,k-1)$ words
$\alpha$ with 0, 1 and 2 copies of $k$, respectively. Hence
\[
f_{111}(n,k)=f_{111}(n,k-1)+nf_{111}(n-1,k-1)+\binom{n}{2}f_{111}(n-2,k-1),
\]
for all $n,k\ge 2$. Also, $f_{111}(n,1)=1$ for $n=0,1,2$,
$f_{111}(n,1)=0$ for all $n\ge 3$, $f_{111}(0,k)=1$ and
$f_{111}(1,k)=k$ for all $k$, hence the theorem holds.
\end{proof}

Finally, we consider patterns 112 and 121. We start with pattern
121.

If a word $\sigma\in[k]^n$ avoids pattern 121, then it contains no
letters other than 1 between any two 1's, which means that all 1's
in $\sigma$, if any, are consecutive. Deletion of all 1's from
$\sigma$ leaves another word $\sigma_1$ which avoids 121 and
contains no 1's, so all 2's in $\sigma_1$, if any, are
consecutive. In general, deletion of all letters 1 through $j$
leaves a (possibly empty) word $\sigma_j$ on letters $j+1$ through
$k$ in which all letters $j+1$, if any, occur consecutively.

If a word $\sigma\in[k]^n$ avoids pattern 112, then only the
leftmost 1 of $\sigma$ may occur before a greater letter. The rest
of the 1's must occur at the end of $\sigma$. In fact, just as in
the previous case, deletion of all letters $1$ through $j$ leaves
a (possibly empty) word $\sigma_j$ on letters $j+1$ through $k$ in
which all occurrences of $j+1$, except possibly the leftmost one,
are at the end of $\sigma_j$. We will call all occurrences of a
letter $j$, except the leftmost $j$, \emph{excess} $j$'s.

The preceding analysis suggests a natural bijection $\rho:
[k]^n(121)\to [k]^n(112)$. Given a word $\sigma\in[k]^n(121)$, we
apply the following algorithm of $k$ steps. Say it yields a word
$\sigma^{(j)}$ after Step $j$, with $\sigma^{(0)}=\sigma$. Then
Step $j$ ($1\le j\le k$) is:

Step $j$. Cut the block of excess $j$'s, then insert it
immediately before the final block of all smaller excess letters
of $\sigma^{(j-1)}$, or at the end of $\sigma^{(j-1)}$ if there
are no smaller excess letters.

It is easy to see that, at the end of the algorithm, we get a word
$\sigma^{(k)}\in[k]^n(112)$.

The inverse map, $\rho^{-1}: [k]^n(112)\to [k]^n(121)$ is given by
a similar algorithm of $k$ steps. Given a word
$\sigma\in[k]^n(112)$ and keeping the same notation as above, Step
$j$ is as follows:

Step $j$. Cut the block of excess $j$'s (which are at the end
of $\sigma^{(j-1)}$), then insert it immediately after the
leftmost $j$ in $\sigma^{(j-1)}$.

Clearly, we get $\sigma^{(k)}\in[k]^n(121)$ at the end of the
algorithm.

Thus, we have the following
\begin{theorem}\label{f112-121b}
Patterns 121 and 112 are Wilf-equivalent.
\end{theorem}

We will now find $f_{112}(n,k)$ and provide a bijective
proof of the resulting formula.

Consider all words $\sigma\in[k]^n(112)$ which contain a letter 1.
Their number is
\begin{equation} \label{g112-1}
g_{112}(n,k)=f_{112}(n,k)-|\{\sigma\in[k]^n(112) : \sigma \mbox{
has no 1's}\}|=f_{112}(n,k)-f_{112}(n,k-1).
\end{equation}

On the other hand, each such $\sigma$ either ends on 1 or not.

If $\sigma$ ends on 1, then deletion of this 1 may produce any
word in $\bar{\sigma}\in[k]^{n-1}(112)$, since addition of the
rightmost 1 to any word in $\bar{\sigma}\in[k]^{n-1}(112)$ does
not produce extra occurrences of pattern 112.

If $\sigma$ does not end on 1, then it has no excess 1's, so its
only 1 is the leftmost 1 which does not occur at end of $\sigma$.
Deletion of this 1 produces a word in
$\bar{\sigma}\in\{2,\ldots,k\}^{n-1}(112)$. Since insertion of a
single 1 into each such $\bar{\sigma}$ does not produce extra
occurrences of pattern 112, for each word
$\bar{\sigma}\in\{2,\ldots,k\}^{n-1}(112)$ we may insert a single
1 in $n-1$ positions (all except the rightmost one) to get a word
$\sigma\in[k]^n(112)$ which contains a single 1 not at the
end.

Thus, we have
\begin{multline} \label{g112-2}
g_{112}(n,k)=f_{112}(n-1,k)+(n-1)|\{\sigma\in[k]^{n-1}(112) : \sigma
\mbox{ has no 1's}\}|=\\
=f_{112}(n-1,k)+(n-1)f_{112}(n-1,k-1).
\end{multline}

Combining (\ref{g112-1}) and (\ref{g112-2}), we get
\begin{equation} \label{f112b-r}
f_{112}(n,k)-f_{112}(n,k-1)=f_{112}(n-1,k)+(n-1)f_{112}(n-1,k-1), \quad n\ge 1,\ k\ge 1.
\end{equation}

The initial values are $f_{112}(n,0)=\delta_{n0}$ for all $n\ge 0$
and $f_{112}(0,k)=1$, $f_{112}(1,k)=k$ for all $k\ge 0$.

Therefore, multiplying (\ref{f112b}) by $y^k$ and summing over
$k$, we get
\[
G_{112}(n;y)-\delta_{n0}-yG_{112}(n;y)=G_{112}(n-1;y)-\delta_{n-1,0}+(n-1)yG_{112}(n-1;y),
\quad n\ge 1,
\]
hence,
\[
(1-y)G_{112}(n;y)=(1+(n-1)y)G_{112}(n-1;y), \quad n\ge 2.
\]
Therefore,
\begin{equation} \label{g112b-r}
G_{112}(n;y)=\frac{1+(n-1)y}{1-y}\,G_{112}(n-1;y), \quad n\ge 2.
\end{equation}
Also, $\displaystyle G_{112}(0;y)=\frac{1}{1-y}$ and
$\displaystyle G_{112}(1;y)=\frac{y}{(1-y)^2}$, so applying the previous
equation repeatedly yields
\begin{equation} \label{g112b}
G_{112}(n;y)=\frac{y(1+y)(1+2y)\cdots(1+(n-1)y)}{(1-y)^{n+1}}.
\end{equation}

We have
\begin{multline*}
\frac{1}{y} Numer(G_{112}(n;y))=(1+y)(1+2y)\cdots(1+(n-1)y)
=y^n\prod_{j=0}^{n-1}{\left(\frac{1}{y}+j\right)}=\\
=y^n\sum_{k=0}^{n}{c(n,k)\left(\frac{1}{y}\right)^k}
=\sum_{k=0}^{n}{c(n,k)y^{n-k}}=\sum_{k=0}^{n}{c(n,n-k)y^k},
\end{multline*}
where $c(n,j)$ is the signless Stirling number of the first kind,
and
\[
y Denom(G_{112}(n;y))=\frac{y}{(1-y)^{n+1}}=
\sum_{k=0}^{\infty}{\binom{n+k-1}{n}y^k},
\]
so $f(n,k)$ is the convolution of the two coefficients:
\[
f_{112}(n,k)=\left(c(n,n-k)\ast\binom{n+k-1}{n}\right)=
\sum_{j=0}^{k}{\binom{n+k-j-1}{n}c(n,n-j)}.
\]
Thus, we have a new and improved version of Theorem
\ref{f112-121b}.
\begin{theorem} \label{f112-121b2}
Patterns 112 and 121 are Wilf-equivalent, and
\begin{equation} \label{f112b}
\begin{split}
f_{121}(n,k)=f_{112}(n,k)&=\sum_{j=0}^{k}{\binom{n+k-j-1}{n}c(n,n-j)},\\
F_{121}(x,y)=F_{112}(x,y)&=\frac{1}{1-y}\cdot \left( \frac{1-y}{1-y-xy} \right)^{1/y}.
\end{split}
\end{equation}
\end{theorem}

We note that this is the first time that Stirling numbers appear
in enumeration of words (or permutations) with forbidden patterns.

\begin{proof}
The first formula is proved above. The second formula can be
obtained as the exponential generating function of $G_{112}(n;y)$
from the recursive equation (\ref{g112b-r}). Alternatively,
multiplying the recursive formula (\ref{f112b-r}) by $x^n/n!$ and
summing over $n$ yields
\[
\frac{d}{dx}F_{112}(x;k)=F_{112}(x;k)+(1+x)\frac{d}{dx}F_{112}(x;k-1).
\]
Multiplying this by $y^k$ and summing over $k\ge 1$, we obtain
\[
\frac{d}{dx} F_{112}(x,y)=\frac{1}{1-y-yx}F_{112}(x,y).
\]
Solving this equation together with the initial condition
$\displaystyle F_{112}(0,y)=\frac{1}{1-y}$ yields the desired
formula.
\end{proof}

We will now prove the exact formula (\ref{f112b}) bijectively. As
it turns out, a little more natural bijective proof of the same
formula obtains for $f_{221}(n,k)$, an equivalent result since
$221=c(112)$. This bijective proof is suggested by equation
(\ref{f112b-r}) and by the fact that $c(n,n-j)$ enumerates
permutations of $n$ letters with $n-j$ right-to-left minima (i.e.
with $j$ right-to-left nonminima), and $\binom{n+k-j-1}{n}$
enumerates nondecreasing strings of length $n$ on letters in
$\{0,1,\ldots,k-j-1\}$.

Given a permutation $\pi\in S_n$ which has $n-j$ right-to-left
minima, we will construct a word $\sigma\in[j+1]^n(221)$ with
certain additional properties to be discussed later. The algorithm
for this construction is as follows.

\begin{algorithm} \label{minword}
~\\[-30pt]
\begin{enumerate}

\item \label{assign1} Let $d=(d_1,\ldots,d_n)$, where
$
d_r=\begin{cases}
0, \text{ if $r$ is a right-to-left minimum in $\pi$,}\\
1, \text{ otherwise.}
\end{cases}
$

\item Let $s=(s_1,s_2,\ldots,s_n)$, where
$s_r=1+\sum_{i=1}^{r}{d_r}$, $r=1,\ldots,n$.

\item Let $\sigma=\pi\circ s$ (i.e. $\sigma_{r}=s_{\pi(r)}$,
$r=1,\ldots,n$). This is the desired word $\sigma$.

\end{enumerate}
\end{algorithm}

\begin{example} \label{exminword}
Let $\pi=621/93/574/8/10\in S_{10}$. Then $n-j=5$, so
$j+1=6$, $d=0100111010$, $s=1222345566$, so the corresponding word
$\sigma=4216235256\in[6]^{10}(221)$.
\end{example}

Note that each letter $s_r$ in $\sigma$ is in the same position as
that of $r$ in $\pi$, i.e. $\pi^{-1}(r)$.

Let us show that our algorithm does indeed produce a word
$\sigma\in[j+1]^n(221)$.

Since $\pi$ has $n-j$ right-to-left minima, only $j$ of the
$d_r$'s are 1s, the rest are 0s. The sequence $\{s_r\}$ is clearly
nondecreasing and its maximum, $s_n=1+1\cdot j=j+1$. Thus,
$\sigma\in[j+1]^n$ and $\sigma$ contains all letters from 1 to
$j+1$.

Suppose now $\sigma$ contains an occurrence of the pattern 221.
This means $\pi$ contains a subsequence $b c a$ or $c b a$,
$a<b<c$. On the other hand, $s_b=s_c$, so
$0=s_c-s_b=\sum_{r=b+1}^{c}{d_r}$, hence $d_c=0$ and $c$ must be a
right-to-left minimum. But $a<c$ is to the right of $c$, so $c$ is
not a right-to-left minimum. Contradiction. Therefore, $\sigma$
avoids pattern 221.

Thus, $\sigma\in[j+1]^n(221)$ and contains all letters 1 through
$j+1$. Moreover, the leftmost (and \emph{only} the leftmost)
occurrence of each letter (except 1) is to the left of some
smaller letter. This is because $s_{b}=s_{b-1}$ means $d_b=0$,
that is $b$ is a right-to-left minimum, i.e. occurs to the right
of all smaller letters. Hence, $s_b$ is also to the right of all
smaller letters, i.e. is a right-to-left minimum of $\sigma$. On
the other hand, $s_b>s_{b-1}$ means $d_b=1$, that is $b$ is not a
right-to-left minimum of $\pi$, so $s_b$ is not a right-to-left
minimum of $\sigma$.

It is easy to construct an inverse of Algorithm 1. Assume we are
given a word $\sigma$ as above. We will construct a permutation
$\pi\in S_n$ which has $n-j$ right-to-left minima.

\begin{algorithm} \label{wordmin}
~\\[-24pt]
\begin{enumerate}

\item Reorder the elements of $\sigma$ in nondecreasing order and
call the resulting string $s$.

\item Let $\pi\in S_n$ be the permutation such that
$\sigma_{r}=s_{\pi(r)}$, $r=1,\ldots,n$, given that
$\sigma_a=\sigma_b$ (i.e. $s_{\pi(a)}=s_{\pi(b)}$) implies
$\pi(a)<\pi(b) \Leftrightarrow a<b$). In other words, $\pi$ is
monotone increasing on positions of equal letters. Then $\pi$ is
the desired permutation.

\end{enumerate}
\end{algorithm}

\begin{example} \label{exwordmin}
Let $\sigma=4216235256\in[6]^{10}(221)$ from our
earlier example (so $j+1=6$). Then $s=1222345566$, so looking at
positions of 1s, 2s, etc., 6s, we get
\begin{eqnarray*}
\pi(1)=6 & & \\
\pi(\{2,5,8\})=\{2,3,4\} &\implies & \pi(2)=2,\ \pi(5)=3,\ \pi(8)=4\\
\pi(3)=1 & & \\
\pi(\{4,10\})=\{9,10\} &\implies & \pi(9)=4,\ \pi(10)=10\\
\pi(6)=5 & & \\
\pi(\{7,9\})=\{7,8\} &\implies & \pi(7)=7,\ \pi(9)=8.
\end{eqnarray*}
Hence, $\pi=(6,2,1,9,3,5,7,4,8,10)$ (in the one-line notation, not
the cycle notation) and $\pi$ has $n-j$ right-to-left minima: 10,
8, 4, 3, 1.
\end{example}

Note that the position of each $s_r$ in $\sigma$ is $\pi^{-1}(r)$,
i.e. again the same as $r$ has in $\pi$. Therefore, we conclude as
above that $\pi$ has $j+1-1=j$ right-to-left nonminima, hence,
$n-j$ right-to-left minima. Furthermore, the same property implies
that Algorithm 2 is the inverse of Algorithm 1.

Note, however, that more than one word in $[k]^n(221)$ may map to
a given permutation $\pi\in S_n$ with exactly $n-j$ right-to-left
minima. We only need require that just the letters corresponding
to the right-to-left nonminima of $\pi$ be to the left of a
smaller letter (i.e. not at the end) in $\sigma$. Values of 0 and
1 of $d_r$ in Step \ref{assign1} of Algorithm 1 are minimal
increases required to recover back the permutation $\pi$ with
Algorithm 2. We must have $d_r\ge 1$ when we \emph{have to}
increase $s_r$, that is when $s_r$ is not a right-to-left minimum
of $\sigma$, i.e. when $r$ is not a right-to-left minimum of
$\pi$. Otherwise, we don't have to increase $s_r$, so $d_r\ge 0$.

Let $\sigma\in[k]^n(221)$, $\pi=Alg2(\sigma)$,
$\tilde\sigma=Alg1(\pi)=Alg1(Alg2(\sigma))\in[j+1]^n(221)$, and
$\eta=\sigma-\tilde\sigma$ (vector subtraction). Note that
$e_r=s_r(\sigma)-s_r(\tilde\sigma)\ge 0$ does not decrease (since
$s_r(\sigma)$ cannot stay the same if $s_r(\tilde\sigma)$ is
increased by 1) and $0\le e_1\le\ldots\le e_n\le k-j-1$.

Since position of each $e_r$ in $\eta$ is the same as position of
$s_r$ in $\sigma$ (i.e. $\eta_a=e_{\pi(a)}$, $e=e_1e_2\ldots
e_n$), the number of such sequences $\eta$ is the number of
nondecreasing sequences $e$ of length $n$ on letters in
$\{0,\ldots,k-j-1\}$, which is $\binom{n+k-j-1}{n}$.

Thus, $\sigma\in[k]^n(221)$ uniquely determines the pair
$(\pi,e)$, and vice versa. This proves the formula (\ref{f112b})
of Theorem \ref{f112-121b2}.

All of the above lets us state the following

\begin{theorem} \label{theorem3}
There are 3 Wilf classes of multipermutations of length 3, with
representatives 123, 112 and 111.
\end{theorem}


\section{Pairs of 3-letter patterns}

There are 8 symmetric classes of pairs of $3$-letters words, which
are
\[
\{111,112\}, \{111, 121\}, \{112,121\}, \{112,122\}, \{112,211\},
\{112,212\}, \{112,221\}, \{121,212\}.
\]

\begin{theorem}
\label{f111112}
The pairs $\{111,112\}$ and $\{111, 121\}$ are Wilf equivalent, and
\[
\begin{split}
F_{111,121}(x,y)=F_{111,112}(x,y)=\frac{e^{-x}}{1-y}\cdot
\left(\frac{1-y}{1-y-xy} \right)^{1/y},\\
f_{111,112}(n,k)=\sum_{i=0}^{n}\sum_{j=0}^{k}{(-1)^{n-i}\binom{n}{i}\binom{k+i-j-1}{i}c(i,i-j)}.
\end{split}
\]
\end{theorem}
\begin{proof}
To prove equivalence, notice that the bijection
$\rho:[k]^n(121)\to[k]^n(112)$ preserves the number of excess
copies of each letter and that avoiding pattern 111 is the same as
having at most 1 excess letter $j$ for each $j=1,\ldots,k$. Thus,
restriction of $\rho$ to words with $\le 1$ excess letter of each
kind yields a bijection
$\rho\!\restriction_{111}:[k]^n(111,121)\to[k]^n(111,112)$.

Let $\alpha\in [k]^n(111,112)$ contain $i$ copies of letter $1$.
Since $\alpha$ avoids $111$, we see that $i\in\{0,1,2\}$.
Corresponding to these three cases, the number of such words
$\alpha$ is $f_{111,112}(n,k-1)$, $nf_{111,112}(n-1,k-1)$ or
$(n-1)f_{111,112}(n-2,k-1)$, respectively. Therefore,
\[
f_{111,112}(n,k)=f_{111,112}(n,k-1)+nf_{111,112}(n-1,k-1)+(n-1)f_{111,112}(n-2,k-1),
\]
for $n,k\ge 1$. Also, $f_{111,112}(n,0)=\delta_{n0}$ and
$f_{111,112}(0,k)=1$, hence
\[
F_{111,112}(x;k)=(1+x)F_{111,112}(x;k-1)+\int\!\!{xF_{111,112}(x;k-1)dx},
\]
where $f_{111,112}(0,k)=1$. Multiply the above equation by $y^k$
and sum over all $k\ge 1$ to get
\[
F_{111,112}(x,y)=c(y)e^{-x}\cdot\left(\frac{1-y}{1-y-xy}\right)^{1/y},
\]
which, together with $\displaystyle
F_{111,112}(0,y)=\frac{1}{1-y}$, yields the generating function.

Notice that $F_{111,112}(x,y)=e^{-x}F_{112}(x,y)$, hence,
$F_{111,112}(x;k)=e^{-x}F_{112}(x;k)$, so $f_{111,112}(n,k)$ is
the exponential convolution of $(-1)^n$ and $f_{112}(n,k)$. This
yields the second formula.
\end{proof}

\begin{theorem}
\label{f112121}
Let $H_{112,121}(x;k)=\sum_{n\ge 0} f_{112,121}(n,k)x^n$. Then for
any $k\ge 1$,
\[
H_k(x)=\frac{1}{1-x}H_{112,121}(x;k-1)+x^2\frac{d}{dx}H_{112,121}(x;k-1),
\]
and $H_{112,121}(x;0)=1$.
\end{theorem}
\begin{proof}
Let $\alpha\in [k]^n(112,121)$ such that contains $j$ letters $1$.
Since $\alpha$ avoids $112$ and $121$, we have that for $j>1$, all
$j$ copies of letter $1$ appear in $\alpha$ in positions $n-j+1$
through $n$. When $j=1$, the single 1 may appear in any position.
Therefore,
\[
f_{112,121}(n;k)=f_{112,121}(n;k-1)+nf_{112,121}(n-1,k-1)+\sum_{j=2}^{n}{f_{112,121}(n-j;k-1)},
\]
which means that
\begin{multline*}
f_{112,121}(n;k)=f_{112,121}(n-1;k)+f_{112,121}(n;k-1)\\+(n-1)f_{112,121}
(n-1,k-1)-(n-2)f_{112,121}(n-2,k-1).
\end{multline*}
We also have $f_{112,121}(n;0)=1$, hence it is easy to see the
theorem holds.
\end{proof}

\begin{theorem}
\label{f112211}
Let $H_{112,211}(x;k)=\sum_{n\ge 0}{f_{112,211}(n,k)x^n}$. Then
for any $k\ge 1$,
\[
H_{112,211}(x;k)=(1+x+x^2)H_{112,211}(x;k-1)+\frac{x^3}{1-x}+\frac{d}{dx}H_{112,211}(x;k-1),
\]
and $H_{112,211}(x;0)=1$.
\end{theorem}
\begin{proof}
Let $\alpha\in [k]^n(112,211)$ such that contains $j$ letters $1$.
Since $\alpha$ avoids $112$ and $211$ we have that $j=0,1,2,n$.
When $j=2$, the two 1's must at the beginning and at the end.
Hence, it is easy to see that for $j=0,1,2,n$ there are
$f_{112,211}(n;k-1)$, $nf_{112,211}(n-1;k-1)$,
$f_{112,211}(n-2;k-1)$ and $1$ such $\alpha$, respectively.
Therefore,
\[
f_{112,211}(n;k)=f_{112,211}(n;k-1)+nf_{112,211}(n-1,k-1)+f_{112,211}(n-2,k-1)+\delta_{n\ge3}.
\]
We also have $f_{112,121}(n;0)=1$, hence it is easy to see the
theorem holds.
\end{proof}

\begin{theorem}
\label{f112212} Let $a_{n,k}=f_{112,212}(n,k)$, then
\[
a_{n,k}=a_{n,k-1}+\sum_{d=1}^n \sum_{r=0}^{k-1} \sum_{j=0}^{n-d}
a_{j,r}a_{n-d-j,k-1-r}
\]
and $a_{0,k}=1$, $a_{n,1}=1$.
\end{theorem}
\begin{proof}
Let $\alpha\in [k]^n(112,212)$ have exactly $d$ letters $1$. If
$d=0$, there are $a_{n,k-1}$ such $\alpha$. Let $d\ge 1$, and
assume that $\alpha_{i_d}=1$ where $d=1,2,\dots j$. Since $\alpha$
avoids $112$, we have $i_{2}=n+2-d$ (if $d=1$, we define
$i_2=n+1$), and since $\alpha$ avoids $212$ we have that
$\alpha_a,\alpha_b$ are different for all $a<i_1<b<i_2$.
Therefore, $\alpha$ avoids $\{112,212\}$ if and only if
$(\alpha_1,\dots,\alpha_{i_1-1})$, and
$(\alpha_{i_1+1},\dots,\alpha_{i_2-1})$ are
$\{112,212\}$-avoiding. The rest is easy to obtain.
\end{proof}

\begin{theorem}
\label{f112221}
\[
f_{112,221}(n,k)=\sum_{j=1}^{k}{j\cdot j!\binom{k}{j}}
\]
for all $n\ge k+1$,
\[
f_{112,221}(n,k)=n!\binom{k}{n}+\sum_{j=1}^{n-1}{j\cdot
j!\binom{k}{j}}
\]
for all $k\ge n\ge 2$, and $f_{112,221}(0,k)=1$,
$f_{112,221}(1,k)=k$.
\end{theorem}
\begin{proof}
Let $\alpha\in [k]^n(112,221)$ and $j\le n$ be such that
$\alpha_1,\ldots,\alpha_j$ are all distinct and $j$ is maximal.
Clearly, $j\le k$. Since $\alpha$ avoids $\{112, 221\}$ and $j$ is
maximal, we get that the letters $\alpha_{j+1},\ldots,\alpha_n$,
if any, must all be the same and equal to one of the letters
$\alpha_1,\ldots,\alpha_j$. Hence, there are $j\cdot
j!\binom{k}{j}$ such $\alpha$ if , for $j<n$ or $j=n>k$. For
$j=n\le k$, there are $n!\binom{k}{n}$ such $\alpha$. Hence,
summing over all possible $j=1,\ldots,k$, we obtain the theorem.
\end{proof}

\begin{theorem}
\label{f121212}
\[
f_{121,212}(n,k)=\sum_{j=0}^{k}{j!\binom{k}{j}\binom{n-1}{j-1}}
\]
for $k\ge 0$, $n\ge 1$, and $f_{121,212}(0,k)=1$ for $k\ge 0$.
\end{theorem}

\begin{proof}
Let $\alpha\in [k]^n(121,212)$ contain exactly $j$ distinct
letters. Then all copies of each letter $1$ through $j$ must be
consecutive, or $\alpha$ would contain an occurrence of either
$121$ or $212$. Hence, $\alpha$ is a concatenation of $j$ constant
strings. Suppose the $i$-th string has length $n_i>0$, then
$n=\sum_{i=1}^{j}{n_i}$. Therefore, to obtain any $\alpha\in
[k]^n(121,212)$, we can choose $j$ letters out of $k$ in
$\binom{k}{j}$ ways, then choose any ordered partition of $n$ into
$j$ parts in $\binom{n-1}{j-1}$ ways, then label each part $n_i$
with a distinct number $l_i\in\{1,\ldots,j\}$ in $j!$ ways, then
substitute $n_i$ copies of letter $l_i$ for the part $n_i$
($i=1,\ldots,j$). This yields the desired formula.
\end{proof}


Unfortunately, the case of the pair $(112,122)$ still remains unsolved.


\section{Some triples of $3$-letter patterns}

\begin{theorem}
\label{f112121211}
\[
\begin{split}
F_{112,121,211}(x;k)&=1+\frac{(e^x-1)((1+x)^k-1)}{x},\\
f_{112,121,211}(n,k)&=
\begin{cases}
\displaystyle\sum_{j=1}^{n}{\frac{1}{j!}\binom{n+1}{j}\binom{k}{n+1-j}}, \quad n\ge 1,\\
1, \quad n=0.
\end{cases}
\end{split}
\]
\end{theorem}
\begin{proof}
Let $\alpha\in [k]^n(112,121,211)$ contain $j$ letters $1$. For
$j\ge 2$, there are no letters between the $1$'s, to the left of
the first $1$ or to the right of the last $1$, hence $j=n$. For
$j=1$, $j=0$  it is easy to see from definition that there are
$nf_{112,121,211}(n-1,k-1)$ and $f_{112,121,211}(n,k-1)$ such
$\alpha$, respectively. Hence,
\[
f_{112,121,211}(n,k)=f_{112,121,211}(n,k-1)+nf_{112,121,211}(n-1,k-1)+1,
\]
for $n,k\ge 2$. Also, $a(n,1)=a(n,0)=1$, $a(0,k)=1$, and
$a(1,k)=k$. Let $b(n,k)=f_{112,121,211}(n,k)/n!$, then
\[
b(n,k)=b(n,k-1)+b(n-1,k-1)+\frac{1}{n!}.
\]
Let $b_k(x)=\sum_{n\ge 0} b(n,k)x^n$, then it is easy to see that
$b_k(x)=(1+x)b_{k-1}(x)+e^x-1$. Since we also have $b_0(x)=e^x$,
the theorem follows by induction.
\end{proof}


\section{Some patterns of arbitrary length}

\subsection{Pattern $11\ldots1$}

Let us denote by $\langle a\rangle_l$ the word consisting of $l$
copies of letter $a$.

\begin{theorem}
For any $l,k\ge 0$,
\[
F_{\langle1\rangle_l}(x;k)=\left( \sum_{j=0}^{l-1} \frac{x^j}{j!}
\right)^k.
\]
\end{theorem}
\begin{proof}
Let $\alpha\in [k]^n(\langle 1\rangle_l)$ contain $j$ letters $1$.
Since $\alpha$ avoids $\langle 1\rangle_l$, we have $j\le l-1$. If
$\alpha$ contains exactly $j$ letters of $1$, then there are
$\binom{n}{j}f_{\langle 1\rangle_l}(n-j,k-1)$ such $\alpha$,
therefore
\[
f_{\langle1\rangle_l}(n,k)=\sum_{j=0}^{l-1}{\binom{n}{j}
f_{\langle1\rangle_l}(n-j,k-1)}.
\]
We also have $f_{\langle 1\rangle_l}(n,k)=k^n$ for $n\le l-1$,
hence it is easy to see the theorem holds.
\end{proof}

In fact, \cite{CS} shows that we have
\[
f_{\langle1\rangle_l}(n,k)=\sum_{i=1}^{n}{M_2^{l-1}(n,i)(k)_i},
\]
where $M_2^{l-1}(n,i)$ is the number of partitions of an $n$-set
into $i$ parts of size $\le l-1$.

\subsection{Pattern $11\ldots 121\ldots 11$}

Let us denote $v_{m,l}=11\dots 121\dots 11$, where $m$
(respectively, $l$) is the number of $1$'s on the left
(respectively, right) side of $2$ in $v_{m,l}$. In this section we
prove the number of words in $[k]^n(v_{m,l})$ is the same as the
number of words in $[k]^n(v_{m+l,0})$ for all $m,l\geq 0$.

\begin{theorem}
\label{mml}
Let $m,l\ge 0$, $k\ge 1$. Then for $n\ge 1$,
\[
f_{v_{m,l}}(n+1,k)-f_{v_{m,l}}(n,k)=\sum_{j=0}^{m+l-1}{\binom{n}{j}f_{v_{m,l}}(n+1-j,k-1)}.
\]
\end{theorem}
\begin{proof}
Let $\alpha\in [k]^n(v_{m,l})$ contain exactly $j$ letters $1$.
Since the $1$'s cannot be part of an occurrence of $v_{m,l}$ in
$\alpha$ when $j\le m+l-1$, these $1$'s can be in any $j$
positions, so there are $\binom{n}{j}f_{v_{m,l}}(n,k-1)$ such
$\alpha$. If $j\ge m+l$, then the $m$-th through $(j-l+1)$-st
($l$-th from the right) $1$'s must be consecutive letters in
$\alpha$ (with the convention that the $0$-th 1 is the beginning
of $\alpha$ and $(j+1)$-st 1 is the end of $\alpha$). Hence, there
are $\binom{n-j+m+l-1}{m+l-1}f_{v_{m,l}}(n-j,k-1)$ such $\alpha$,
and hence
\[
f_{v_{m,l}}(n;k)=\sum_{j=0}^{m+l-1}{\binom{n}{j}{f_{v_{m,l}}(n-j,k-1)}}
 + \sum_{j=m+l}^{n}{\binom{n-j+m+l-1}{m+l-1}f_{v_{m,l}}(n-j,k-1)}.
\]
Hence for all $n\ge 1$,
\[
f_{v_{m,l}}(n+1,k)-f_{v_{m,l}}(n,k)=\sum_{j=0}^{m+l-1}{\binom{n}{j}f_{v_{m,l}}(n+1-j,k-1)}.
\]
\end{proof}

An immediate corollary of Theorem \ref{mml} is the following.

\begin{corollary} \label{vml}
Let $m,l\ge 0$, $k\ge 0$. Then for $n\ge 0$
\[
f_{v_{m,l}}(n,k)=f_{v_{m+l,0}}(n,k).
\]
In other words, all patterns $v_{m,l}$ with the same $m+l$ are
Wilf-equivalent.
\end{corollary}

\begin{proof}
We will give an alternative, bijective proof of this by
generalizing our earlier bijection $\rho:[k]^n(121)\to[k]^n(112)$.
Let $\alpha\in [k]^n(v_{m,l})$. Recall that $\alpha_j$ is a word
obtained by deleting all letters $1$ through $j$ from $\alpha$
(with $\alpha_0:=\alpha$).

Suppose that $\alpha$ contains $i$ letters $j+1$. Then all
occurrences of $j+1$ from $m$-th through $(i-l+1)$-st, if any
(i.e. if $j\ge m+l$), must be consecutive letters in $\alpha_j$.
We will denote as \emph{excess} $j$'s the $(m+1)$-st through
$(i-l+1)$-st copies of $j$ when $l>0$, and $m$-th through $i$-th
copies of $j$ when $l=0$.

Suppose that $m+l=m'+l'$. Then the bijection
$\rho_{m,l;m'\!,l'}:[k]^n(v_{m,l})\to[k]^n(v_{m'\!,l'})$ is an
algorithm of $k$ steps. Given a word $\alpha\in[k]^n(v_{m,l})$,
say it yields a word $\alpha^{(j)}$ after Step $j$, with
$\alpha^{(0)}:=\alpha$. Then Step $j$ ($1\le j\le k$) is as
follows:

Step $j$.\\[-24pt]
\begin{enumerate}
\item Cut the block of excess $j$'s from ${\alpha^{(j-1)}}_{j-1}$ (which is
immediately after the $m$-th occurrence of $j$), then insert it
immediately after the $m'$-th occurrence of $j$ if $l'>0$, or at
the end of ${\alpha^{(j-1)}}_{j-1}$ if $l'=0$.

\item Insert letters $1$ through $j-1$ into the resulting string
in the same positions they are in $\alpha^{(j-1)}$ and call the
combined string $\alpha^{(j)}$.
\end{enumerate}

Clearly,
\[
{\alpha^{(j)}}_{j}={\alpha^{(j-1)}}_{j}=\ldots={\alpha^{(0)}}_{j}=\alpha_j
\]
and at Step $j$, the $j$'s are rearranged so that no $j$ can be
part of an occurrence of $v_{m'\!,l'}$. Also, positions of letters
$1$ through $j-1$ are the same in $\alpha^{(j)}$ and
$\alpha^{(j-1)}$, hence, no letter from $1$ to $j$ can be part of
$v_{m'\!,l'}$ in $\alpha^{(j)}$ by induction. Therefore,
$\alpha^{(k)}\in[k]^n(v_{m'\!,l'})$ as desired.

Clearly, this map is invertible, and
$\rho_{m'\!,l';m,l}=(\rho_{m,l;m'\!,l'})^{-1}$. This ends the
proof.
\end{proof}


\begin{theorem}
\label{thp}
Let $p\ge 1$ and $d_p(f(x))=\int\dots\int f(x)dx\dots
dx$ (and we define $d_0(f(x))=f(x))$). Then for any $k\ge 1$,
\[
F_{v_{p,0}}(x;k)-\int F_{v_{p,0}}(x;k)dx = \sum_{j=0}^{p-1}
{\left( (-1)^jd_p(F_{v_{p,0}}(x;k-1)) \sum_{i=0}^{p-1-j}
{\frac{x^i}{i!}} \right)},
\]
and $F_{v_{p,0}}(x;1)=e^x$, $F_{v_{p,0}}(0;k)=1$.
\end{theorem}

\begin{proof}
By definition, we have $f_{v_{p,0}}(n,1)=1$ for all $n\ge 0$ so
$F_{v_{p,0}}(x;1)=e^x$. On the other hand, Theorem \ref{mml}
yields immediately the rest of this theorem.
\end{proof}

\begin{example}
For $p=1$, Theorem \ref{thp} yields
\[
\sum_{n\ge 0} |[k]^n(12)|\frac{x^n}{n!}=e^x\sum_{j=0}^{k-1}{\binom{k-1}{j}\frac{x^j}{j!}},
\]
which means that, for any $n\ge 0$
\[
|[k]^n(12)|=\binom{n+k-1}{k-1}.
\]
(cf. Section 2.)
\end{example}

\begin{example}
For $p=2$, Theorem \ref{thp} yields
\[
F_{112}(x;k)=e^x\cdot \int (1+x)e^{-x}F_{112}(x;k-1)dx,
\]
and $F_{112}(x;0)=1$.
\end{example}

\begin{corollary}
For any $p\geq 0$
\[
F_{v_{p,0}}(x;2)=e^x\sum_{j=0}^p \frac{x^j}{j!}.
\]
\end{corollary}
\begin{proof}
From Theorem \ref{thp}, we immediately get that
\[
F_{v_{p,0}}(x;2)-\int F_{v_{p,0}}(x;2)dx=e^x\sum_{j=0}^{p-1}{(-1)^j\sum_{i=0}^{p-1-j}{\frac{x^i}{i!}}},
\]
which means that
\[
e^x\frac{d}{dx} \left( e^{-x}F_{v_{p,0}}(x;2) \right)=e^x\sum_{j=0}^{p-1} \frac{x^j}{j!},
\]
hence the corollary holds.
\end{proof}


\end{document}